\magnification=\magstep 1
\baselineskip=15pt     
\parskip=3pt plus1pt minus.5pt
\overfullrule=0pt
\font\hd=cmbx10 scaled\magstep1
\font\hd=cmbx10 scaled\magstep1
\input cyracc.def
\newfam\cyrfam

\def\num{\global\advance\count10 by 1 \eqno(\the\count10)}

\def\Pone{{\bf P}^1}
\def\P{{\bf P}}

\def\X{{\cal X}}
\def\O{{\cal O}}

\def\T{{\cal T}}
\def\H{{\cal H}}

\def\U{{\cal U}}

\def\E{{\cal E}}

\def\EY{\E(1)|_Y}

\def\PEY{\P(\EY)}

\def\OP{\O_{\PEY}(1)}
\def\OP1{\O_{\Pone}}

\def\Pic{{\rm Pic}}
\def\Spec{\mathop{\rm Spec}}

\def\im{\mathop{\rm Im}}

\def\mod{\mathop{\rm mod}}

\def\hom{{\rm Hom}}

\def\boldz{{\bf Z}}

\def\dual#1{{#1}^{\scriptscriptstyle \vee}}

\def\mapright#1{\smash{
  \mathop{\longrightarrow}\limits^{#1}}}

\def\mapdown#1{\Big\downarrow
   \rlap{$\vcenter{\hbox{$\scriptstyle#1$}}$}}

\input epsf.tex
\input amssym.tex

\newsymbol\SEMI226F
\def\rtimes{\mathop{\SEMI}}
\centerline{\hd Examples of Special Lagrangian Fibrations.}
\medskip
\centerline{\it Mark Gross\footnote{*}{Supported in part by the NSF 
and EPSRC.}}
\medskip
\centerline{November 30th, 2000}
\medskip
\centerline{Mathematics Institute}
\centerline{University of Warwick}
\centerline{Coventry, CV4 7AL}
\centerline{mgross@maths.warwick.ac.uk}
\bigskip
\bigskip
{\hd \S 0. Introduction.}

Of late there has been a great deal of interest in special Lagrangian
submanifolds and manifolds fibred in special Lagrangian submanifolds,
motivated by the Strominger-Yau-Zaslow conjecture [33]. One of the
basic approaches to finding examples is to exploit symmetries of the
ambient manifold. If $X$ is a (non-compact) $n$-dimensional Calabi-Yau
manifold with K\"ahler form $\omega$ and nowhere vanishing holomorphic
$n$-form $\Omega$, and if there is an action of a Lie group $G$ on
$X$ preserving these two forms, one can look for $G$-invariant
special Lagrangian submanifolds of $X$. For us, $G$ will be a torus
$T^m$. Using this sort of symmetry to search for examples reduces the
special Lagrangian equations to simpler ones which can be solved.

This technique has been used independently by M. Haskins and D. Joyce
([13] and [21])
to find new examples of special Lagrangian cones and submanifolds, 
while it has been
used independently by E. Goldstein and myself to construct examples
of special Lagrangian fibrations on non-compact Calabi-Yau manifolds
with $T^{n-1}$ actions. This paper is an extended version of an
informally distributed preprint (essentially just the first two sections
of this paper) released at the same time as Goldstein's preprint
[7]. Goldstein has developed his examples in some very interesting
directions somewhat orthogonal to the ones taken here. While
there is some overlap between the examples considered here (especially
in \S 2) and those in [7], our goal will be to develop more global
information about these fibrations.

The first set of examples, which we discuss in \S 2,
are special Lagrangian fibrations on crepant
resolutions of toric Gorenstein singularities. Such examples were already
mentioned in [7]. However, earlier, in [10], I gave
topological fibrations on such crepant resolutions, with the belief
that these would resemble the actual special Lagrangian fibrations
on these manifolds. What is new here is that we show 
this is indeed the case; with some mild
hypotheses on the K\"ahler metric $\omega$ (which we do not require
to be Ricci-flat) we find that the topological construction of \S 3
of [10] coincides with the special Lagrangian fibrations we give here.
This gives us a global understanding of the structure of these fibrations,
complementing the results in [7]. We also give a new variant of this
construction which yields proper fibrations.

The other main set of examples, which we consider in \S 3,
are smoothings (by flat
deformations) of isolated toric Gorenstein singularities. The geometry
of such smoothings are controlled by combinatorics of the toric data,
by results of K. Altmann [1]. These are new examples of special Lagrangian
fibrations not discussed elsewhere
in the literature. It is great fun to see how special Lagrangian fibrations
change if one starts with a crepant resolution of a toric singularity,
degenerates by contracting down to the toric singularity, and then smooths
(this process is often called an extremal transition). We discuss some
examples of this in \S 3.

\S\S 4 and 5 are more speculative in nature. In \S 4, we amplify a
brief discussion from [10] about the connections between special 
Lagrangian fibrations on crepant resolutions of toric singularities
and the local mirror symmetry of [4]. In doing so, we make a connection
with the work of W.-D. Ruan, who came to the description of torus
fibrations on Calabi-Yau hypersurfaces in
toric varieties via the dual picture to the fibrations
on crepant resolutions. We discuss how one might use an $S^1$
symmetry to construct special Lagrangian fibrations dual to the
fibrations on crepant resolutions. It is likely, however, that this
construction will come up against a phenomenon revealed in a very
recent preprint of Joyce [22], showing that in the $S^1$-invariant case,
we might expect to have codimension one rather than codimension two
discriminant loci. This is a serious issue for the SYZ conjecture.
In the last section, we propose a weaker version of the SYZ
conjecture which will hopefully sidestep these issues.

\medskip
{\it Acknowledgments.} I would like to thank N. Hitchin, D. Joyce,
D. Morrison, and P.M.H. Wilson for useful conversations.

\bigskip
{\hd \S 1. The Basic Construction.}

Recall from [12]:

\proclaim Definition 1.1. Let $X$ be a complex $n$-dimensional manifold,
with a Hermitian metric with K\"ahler form $\omega$ and nowhere
vanishing holomorphic $n$-form $\Omega$. Then we say $M\subseteq
X$ is {\it special Lagrangian with respect to $\omega,\Omega$} if $\dim_{\bf R}
M=n$ and $\omega|_M=0$, $\im\Omega|_M=0$.

Note we do not assume that either $d\omega=0$ or the metric is Ricci-flat.
However, if the volume form $\omega^n/n!$ is proportional to
$\Omega\wedge\bar\Omega$, then special Lagrangian submanifolds are
volume minimizing as remarked in [8]; see also [6]. Following Joyce,
if $d\omega=0$, we
will call the triple $(X,\omega,\Omega)$ an {\it
almost Calabi-Yau manifold}. Also, we typically want
to allow singularities in $M$; technically this should be done in the language
of currents, but we won't worry about such technicalities. Instead, just view
$M$ as being a closed set which is a manifold on an open dense subset of $M$,
and the special Lagrangian condition then is required to hold where $M$ is
a manifold.

The main tool we will use for constructing examples of special Lagrangian
fibrations will be the following
result.  This theorem first appeared in print in [7], and similar
results appear in [21]. For completeness, we also give
the proof here.

\proclaim Theorem 1.2. Let $(X,\omega,\Omega)$ be an almost Calabi-Yau
manifold,
and suppose there is an effective action of $T:=T^m$ on $X$
preserving $\omega$ and $\Omega$. Let $\mu_0:X\rightarrow {\bf t}^*=
{\bf R}^m$ be the moment map associated to this action. 
\item{(1)} Let $X_1,\ldots,X_m$ be a basis for the
vector fields generating the action of $T$. Then $\Omega_{red}
=\iota(X_1,\ldots,X_m)\Omega$ descends to an $n-m$-form 
on non-singular points of $Z_p:=\mu_0^{-1}(p)/T$ for $p\in\mu_0(X)$.
\item{(2)} Let $\omega_{red}$ be the induced symplectic form on $Z_p$
via symplectic reduction. If $M_{red}\subseteq Z_p$, let $M$ denote
the pull-back of $M_{red}$ to $\mu_0^{-1}(p)\subseteq X$. Suppose
$M$ is not contained in the set of critical points of $\mu_0$. Then if
$M_{red}\subseteq Z_p$ is special Lagrangian with respect to 
$\omega_{red},\Omega_{red}$, 
$M\subseteq X$ is special Lagrangian with respect to $\omega,\Omega$.
\item{(3)} Suppose that $g:X\rightarrow Y$
is a continuous
map to an $n-m$-dimensional real manifold $Y$, satisfying
$g(t\cdot x)=g(x)$ for $t\in T$. Then the map $f=(\mu_0,g):X
\rightarrow {\bf t}^*\times Y$ has special Lagrangian fibres 
with respect to $\omega,\Omega$ if the induced maps $g:Z_p\rightarrow Y$
have special Lagrangian fibres with respect to $\omega_{red},\Omega_{red}$
for $p$ in a dense subset of $\mu_0(X)$.

Proof. Let $x\in\mu_0^{-1}(p)$, and suppose $x$ is not a critical
point for $\mu_0$. Then the tangent space $T_x Z_p$ of $Z_p$ at the point
represented by $x$ is identified with
$${T_x\mu_0^{-1}(p)\over 
(T_x\mu_0^{-1}(p))^{\omega}}.$$
Now 
$(T_x\mu_0^{-1}(p))^{\omega}$ is the tangent space to the orbit $T\cdot x$,
and this is generated by the tangent vectors $X_1,\ldots,X_m$. 
Thus $\Omega_{red}$ vanishes on 
$(T_x\mu_0^{-1}(p))^{\omega}$. Since $\Omega_{red}$ is invariant under
the action of $T$,  $\Omega_{red}$ descends to an $n-m$-form on
$T_x Z_p$, proving (1). 

Now let $M_{red}\subseteq Z_p$ be special Lagrangian with respect to
$\omega_{red},\Omega_{red}$ with $M$ not contained in the set of critical points
of $\mu_0$. Let $x\in M$ be a regular point for $\mu_0$. Let $Y_1,\ldots,
Y_{n-m}$ be lifts of a basis of tangent vectors to $M_{red}$
at $x\mod T$. Then $T_x M$ has a basis $X_1,\ldots,X_m,Y_1,\ldots,Y_{n-m}$.
By assumption $\omega(X_i,X_j)=0$, and $\omega(Y_i,Y_j)=\omega_{red}(Y_i,Y_j)
=0$ since $M_{red}$ is Lagrangian. Finally, since $Y_i\in T_x\mu_0^{-1}(p)$,
$\omega(X_i,Y_j)=0$ since $X_i\in (T_x\mu_0^{-1}(p))^{\omega}$. Thus
$\omega|_{T_xM}=0$.

To show the fibre is special Lagrangian at $x$, we just observe that
$$\im\Omega(X_1,\ldots,X_m,Y_1,\ldots,Y_{n-m})
=\im\Omega_{red}(Y_1,\ldots,Y_{n-m})=0.$$
Thus shows $M$ is special Lagrangian. Item (3) now follows immediately from
(2).
$\bullet$ 

\bigskip
As a basic application of this,
let $\varphi(x_1,\ldots,x_n)$ be a real-valued function
on an open subset of ${\bf R}^n$ such that $\varphi(|z_1|^2,\ldots,|z_n|^2)$
is pluri-subharmonic on an open subset $U$ of ${\bf C}^n$. Let
$\varphi_i=\partial\varphi/\partial x_i$, $\varphi_{ij}=\partial^2\varphi
/\partial x_i\partial x_j$. Then the induced K\"ahler form is
$$\eqalign{\omega&=
{i\over 2} \partial\bar\partial\varphi(|z_1|^2,\ldots,|z_n|^2)\cr
&={i\over 2} \partial\left(\sum_i\varphi_iz_id\bar z_i\right)\cr
&={i\over 2}\left( \sum_i\varphi_i dz_i\wedge d\bar z_i
+\sum_{i,j} \varphi_{ij} \bar z_j z_i dz_j\wedge d \bar z_i\right).\cr}$$

\proclaim Corollary 1.3. The fibres of the map $f:U\rightarrow {\bf R}^n$
given by
$$f=(\varphi_1|z_1|^2-\varphi_2|z_2|^2,\ldots,\varphi_1 |z_1|^2-\varphi_n
|z_n|^2,\im (i^{n+1}\prod_j z_j))$$
are special Lagrangian with respect to $\omega$ and
$\Omega=dz_1\wedge\cdots\wedge dz_n$.

Proof. $\omega$ is invariant under the natural $T^n$ action on ${\bf C}^n$. This
$T^n$ action is induced by the vector fields
$$X_j=2i\left(\bar z_j{\partial\over\partial \bar z_j}-
z_j{\partial\over\partial z_j}\right),\quad j=1,\ldots,n.$$
Now
$$\eqalign{\iota(X_j)\omega =&\varphi_j(z_jd\bar z_j+\bar z_j dz_j)\cr
&+\sum_i (\varphi_{ij}|z_j|^2z_id\bar z_i+\varphi_{ij} |z_j|^2\bar z_idz_i)\cr
=& d(\varphi_j|z_j|^2).\cr}$$
Thus $X_j$ is a Hamiltonian vector field with Hamiltonian $\varphi_j |z_j|^2$.
In particular, $X_1-X_2,\ldots,X_1-X_n$ generate a $T^{n-1}$ action
with moment map 
$$\mu_0=(f_1,\ldots,f_{n-1})=(\varphi_1|z_1|^2-\varphi_2|z_2|^2,\ldots,
\varphi_1|z_1|^2-\varphi_n|z_2|^2).$$
Furthermore, $f_n=\im(i^{n+1} \prod z_j)$ is constant on the Hamiltonian
trajectories of the first $n-1$ functions. We can now apply
Theorem 1.2 with $\mu_0=(f_1,\ldots,f_{n-1})$ and $g=f_n$, 
as the $T^{n-1}$-action
preserves $\Omega$ also.
Now
$$\eqalign{\iota(X_{f_1},\ldots,X_{f_{n-1}})\Omega
&=(-2i)^{n-1}\iota(z_1{\partial\over\partial z_1}-z_2
{\partial\over\partial z_2},\ldots,z_1{\partial\over\partial z_1}
-z_n{\partial\over\partial z_n})dz_1\wedge\cdots\wedge dz_n\cr
&=\pm(-2i)^{n-1}d(z_1\cdots z_n).\cr}$$
It is then clear that $g$ induces a special Lagrangian fibration on
the surface $\mu_0^{-1}(p)/T^{n-1}$ for all $p$.
$\bullet$
\bigskip

{\hd \S 2. Resolutions of Toric Singularities.}

Let $N\cong \boldz^n$, and let $M=\hom(N,\boldz)$ be the dual
lattice. Put $N_{\bf R}:=N\otimes_{\boldz}{\bf R}$,
$T_{\bf C}(N)=N\otimes_{\boldz}{\bf C}^*$, $T(N)=N\otimes_{\boldz}{\bf R}
/N$. Then $M$ can be naturally identified with the group of characters
${\rm Hom}(T_{\bf C}(N),{\bf C}^*)$, and we will often identify elements of
$M$ with such functions.

Let $\sigma\subseteq
N_{\bf R}$ be a strongly convex rational polyhedral cone. We will
assume that $\sigma$ is a Gorenstein canonical cone. This means that
if $n_1,\ldots,n_s\in N$ are the set of generators of 1-dimensional
faces of $\sigma$, then there exists an $m_0\in M$ such that
$\langle m_0,n_i\rangle=1$ for all $i$ and $\langle m_0,n\rangle\ge 1$
for all $n\in\sigma\cap (N-\{0\})$. Denote by $Y_{\sigma}$ the
corresponding affine toric variety. $Y_{\sigma}$ has Gorenstein canonical
singularities. Let $P$ be the convex hull of $n_1,\ldots,n_s$ in the 
hyperplane $\langle m_0,\cdot\rangle=1$. From now on we will assume
there is a triangulation of $P$ such that the fan $\Sigma$ obtained
as the cone over this triangulation yields a non-singular toric
variety $Y_{\Sigma}$. Then $Y_{\Sigma}\rightarrow Y_{\sigma}$ is a crepant
resolution, and $K_{Y_{\Sigma}}=0$.

In [10], we constructed a topological fibration on $Y_{\Sigma}$. We recall
the construction here. Note that $T(N)$ acts naturally on $Y_{\Sigma}$.
If $N_{m_0}=\{n\in N|\langle m_0,n\rangle=0\}$, then the subtorus
$T(N_{m_0})$ of $T(N)$ also acts on $Y_{\Sigma}$. Then one chooses a
commutative diagram
$$\matrix{Y_{\Sigma}&\mapright{=}&Y_{\Sigma}\cr
\mapdown{q_1}&&\mapdown{\alpha_1\circ q_1}\cr
Y_{\Sigma}/T(N_{m_0})&\mapright{\alpha_1}&{\bf C}\times {\bf R}^{n-1}\cr
\mapdown{q_2}&&\mapdown{(z,x)\mapsto (|z|,x)}\cr
Y_{\Sigma}/T(N)&\mapright{\alpha_2}&{\bf R}_{\ge 0}\times{\bf R}^{n-1}\cr}
\leqno{(2.1)}$$
of homeomorphisms $\alpha_1,\alpha_2$, with $q_1,q_2$ the quotient maps.
If one composes $\alpha_1\circ q_1$ with the map $(z,x)\mapsto (\im z,
x)$, one obtains a map $f:Y_{\Sigma}\rightarrow {\bf R}^n$. This is a 
topological $T^{n-1}\times{\bf R}$ fibration, and its discriminant
locus was analyzed in [10].

Choose a basis $e_1,\ldots,e_n$ of $N$ with dual basis $e_1^*,\ldots,e_n^*$
such that $m_0=e_1^*+\cdots+e_n^*$.
The dual basis corresponds to coordinates $z_1,\ldots,z_n$ on $T_{\bf C}(N)
=N\otimes_{\boldz} {\bf C}^*\subseteq Y_{\Sigma}$.

\proclaim Proposition 2.1. If $\Omega$ is a nowhere vanishing holomorphic
$n$-form on $Y_{\Sigma}$, and $T_{\bf C}(N)$ is identified 
with the unique dense orbit of $T_{\bf C}(N)$ acting on $Y_{\Sigma}$,
then $\Omega|_{T_{\bf C}(N)}= C dz_1\wedge\cdots\wedge dz_n$, where
$C\in {\bf C}$ is a constant.

Proof. The proof is standard: see [28]. We give the complete
proof here. Clearly $\Omega|_{T_{\bf C}(N)}=f dz_1\wedge\cdots\wedge dz_n$
for some holomorphic function $f$. Furthermore, to guarantee that
$\Omega$ has no zeroes on $T_{\bf C}(N)$, $f$ must be a monomial, i.e.
a constant times a character. It is also easy to check that the expression
$dz_1\wedge\cdots\wedge dz_n$ is independent of the choice of basis
$e_1,\ldots,e_n$ subject to the constraint that $m_0=e_1^*+\cdots+e_n^*$.
So we can take $e_1,\ldots,e_n$ to be chosen to be edges of an $n$-dimensional
cone $\tau$ in $\Sigma$; since $Y_{\Sigma}$ is smooth, this forms
a basis, and since each $e_i$ then satisfies $\langle m_0,e_i\rangle =1$,
$m_0=e_1^*+\cdots+e_n^*$. Then $\dual{\tau}$ is generated by
$e_1^*,\ldots,e_n^*$, so $Y_{\Sigma}$ contains an open affine subset $\Spec
{\bf C}[\dual{\tau}\cap M]=\Spec {\bf C}[z_1,\ldots,z_n]$. Then 
$fdz_1\wedge\cdots\wedge dz_n$ extends to a non-zero $n$-form on
this open affine subset if and only if $f$ is constant, as desired.
$\bullet$

\proclaim Theorem 2.2. Let $\omega$ be the K\"ahler form of a K\"ahler
metric on $Y_{\Sigma}$, invariant under the action
of $T(N)$. Let $\Omega$ be the nowhere vanishing $n$-form on $Y_{\Sigma}$
which restricts to $dz_1\wedge\cdots\wedge dz_n$ on $T_{\bf C}(N)$.
Let $\mu:Y_{\Sigma}\rightarrow {\bf R}^n$ be the moment map
associated to this $T(N)$ action, and let $\mu_0:Y_{\Sigma}\rightarrow 
{\bf R}^{n-1}$ be the moment map associated to the $T(N_{m_0})$ action.
Then the function $g:T_{\bf C}(N)\rightarrow {\bf R}$ given by
$$g(z_1,\ldots,z_n)=\im i^{n+1}\prod z_i$$
extends to a map $g:Y_{\Sigma}\rightarrow {\bf R}$. Furthermore,
$f=(g,\mu_0):Y_{\Sigma}\rightarrow {\bf R}^n$ is a special Lagrangian
fibration with respect to $\omega,\Omega$. 
If $\mu$ is proper, then $f$ coincides topologically with
the construction given in [10], \S 3.

Proof. 
We can identify $m_0$ with the character $\prod z_i$ on $T_{\bf C}(N)$.
A priori $m_0$ extends to only a rational function on $Y_{\Sigma}$. In fact,
it extends to a regular function. To show this, we need to show it
extends across every prime divisor of $Y_{\Sigma}$ contained in
$Y_{\Sigma}-T_{\bf C}(N)$. Let $n$ generate a ray in the fan $\Sigma$,
corresponding to some such divisor $D_n$. The dual cone to
$\tau={\bf R}_{\ge 0}n$ is the half-plane
$$\dual{\tau}=\{m\in M| \langle m,n\rangle \ge 0\}.$$
Now $\langle m_0,n\rangle =1$, so $m_0\in \dual{\tau}, -m_0\not\in
\dual{\tau}$. The open affine piece of $Y_{\Sigma}$ corresponding
to the cone $\tau$ is $\Spec {\bf C}[\dual{\tau}\cap M]$. Thus $m_0$
is a regular function on this open set. Since this open set contains
a dense subset of $D_n$, $m_0$ extends across $D_n$,
and in fact takes the value zero on $D_n$, since $-m_0$ (corresponding
to the character $\prod z_j^{-1}$) is not in $\dual{\tau}$. 

Thus $i^{n+1}m_0$ gives a map $Y_{\Sigma}\rightarrow {\bf C}$, from which the
first claim follows. Using the
moment map $\mu_0$ of the $T(N_{m_0})$-action,
we obtain a map $(i^{n+1}m_0,\mu_0):
Y_{\Sigma}\rightarrow {\bf C}\times {\bf R}^{n-1}$. 
Composing this map with $(z,x)\mapsto (\im z,x)$, we obtain
$f:Y_{\Sigma}\rightarrow {\bf R}^n$, an extension of the special Lagrangian
fibration $f:T_{\bf C}(N) \rightarrow {\bf R}^n$ arising in Corollary
1.3.

To finish, we show $f$ coincides with the construction of [10] if $\mu$ is 
proper.

Because we are assuming $\mu$ is proper, the following facts follow from
[15]
Theorem 4.1: $\mu(Y_{\Sigma})$ is convex, and $\mu$ has connected
fibres. Thus, in particular, $\mu$ identifies $\mu(Y_{\Sigma})$ with
$Y_{\Sigma}/T(N)$. Furthermore, $\mu(Y_{\Sigma})$ is a closed,
locally polyhedral convex set, and the extremal points of
$\mu(Y_{\Sigma})$ are images of fixed points of the $T(N)$-action.
Finally the tangent ``wedge'' to $\mu(Y_{\Sigma})$ at such a point
$\mu(x)$ is generated by the weights of the $T(N)$-action
on $T_x Y_{\Sigma}$. There is a 1-1 correspondence between fixed points
of the $T(N)$-action on $Y_{\Sigma}$ and maximal cones $\tau$ of the fan
$\Sigma$. Now since $Y_{\Sigma}$ is non-singular, each such cone
$\tau$ is generated by a basis $v_1,\ldots,v_n$ of $N$, and the weights
of the $T(N)$-representation on $T_x Y_{\Sigma}$, $x$ the point corresponding
to $\tau$, are $v_1^*,\ldots,v_n^*$. Furthermore, since $\langle m_0,v_i
\rangle=1$, it follows that $m_0=\sum v_i^*$. In particular, $m_0$ is
in the interior of the tangent ``wedge'' of each such extremal point.

Let $r:{\bf R}^n\rightarrow {\bf R}^{n-1}$ be given by
$(x_1,\ldots,x_n)\mapsto (x_1-x_2,\ldots,x_1-x_n)$.
The composition $r\circ\mu:Y_{\Sigma}\rightarrow {\bf R}^{n-1}$
is the moment map $\mu_0$ of the $T(N_{m_0})$ action on $Y_{\Sigma}$.

We will now show that if $L=r^{-1}(c)$, $c\in {\bf R}^{n-1}$,
is a line in ${\bf R}^n=M\otimes_{\boldz}
{\bf R}$ parallel to $m_0$, then $L\cap \mu(Y_{\Sigma})$ is a ray.
Indeed, if $L\cap \partial\mu(Y_{\Sigma})$ is non-empty, the 
description above of the tangent wedges of the extremal points
shows that $L\cap\partial\mu(Y_{\Sigma})$ consists of one point.
Thus $L\cap\mu(Y_{\Sigma})$ is closed (as $\mu$ is proper) and has
one boundary point, so it is a ray. If $L\cap \partial\mu(Y_{\Sigma})$
is empty, then either $L\cap\mu(Y_{\Sigma})=\phi$ or $L\subseteq \mu
(Y_{\Sigma})$. In either case, choose a line $l\subseteq {\bf R}^{n-1}$
such that the plane $r^{-1}(l)$ contains $L$, $r^{-1}(l)\cap\mu(Y_{\Sigma})
\not=\phi$, and $r^{-1}(l)\not\subseteq \mu(Y_{\Sigma})$. Then 
$S=r^{-1}(l)\cap\mu(Y_{\Sigma})$ is a closed convex set. It can only
contain $L$ if $S$ is a half-plane with edge parallel to $L$, 
contradicting the description of the tangent wedges to $\mu(Y_{\Sigma})$.
If $L\cap S=\phi$, then there is a supporting line to $S$, $L'$,
parallel to $L$. This also contradicts the description of the tangent
wedges. Thus $L\cap \mu(Y_{\Sigma})$ is a ray.

We now replace $r$ with its restriction to $\mu(Y_{\Sigma})$. We have
$r:\mu(Y_{\Sigma})\rightarrow {\bf R}^{n-1}$
is surjective, with each fibre being a ray parallel to $m_0$.

We now wish to define homeomorphisms $\alpha_1,\alpha_2$ as in the diagram
(2.1) so that 
$f=\im\circ\alpha_1\circ q_1$, where $\im:{\bf C}\times
{\bf R}^{n-1}\rightarrow {\bf R}^n$ is the map $(z,x)\mapsto (\im z,x)$. 
This will show
that $f$ coincides with the construction of [10], \S 3.

Note that the value of $|m_0|$ only depends on the $T(N)$-orbit
of $(z_1,\ldots,z_n)$, so $|m_0|$ descends to a map $|m_0|:\mu(Y_{\Sigma})
\rightarrow {\bf R}_{\ge 0}$. We define $\alpha_2$ as the product map
$(|m_0|,r):\mu(Y_{\Sigma})\rightarrow {\bf R}_{\ge 0}\times {\bf R}^{n-1}$.
To show that $\alpha_2$ is a homeomorphism, it is enough to show that
$|m_0|:r^{-1}(c)\rightarrow {\bf R}_{\ge 0}$ for $c\in {\bf R}^{n-1}$
is a homeomorphism. Now $r^{-1}(c)$ is a ray whose endpoint is in
$\partial\mu(Y_{\Sigma})$, and $\mu^{-1}(\partial\mu(Y_{\Sigma}))=Y_{\Sigma}
-T_{\bf C}(N)$. Thus $|m_0|$ takes the value zero on the endpoint of
$r^{-1}(c)$. So all we need to know is that $|m_0|$ is monotonically
increasing on $r^{-1}(c)$, and goes to $\infty$. 

To see this, note that $Z_c:=\mu_0^{-1}(c)/T(N_0)$ is in fact isomorphic
to ${\bf C}$, and $m_0$ descends to a holomorphic function on $Z_c$.
The maximum modulus theorem then tells us that $|m_0|$ is monotonically
increasing on $r^{-1}(c)$, 
and Liouville's theorem tells us $|m_0|$ goes to $\infty$.

Thus $\alpha_2$ is an homeomorphism.

Similarly, since $(i^{n+1}m_0,\mu_0)$ 
is constant on $T(N_{m_0})$-orbits, it descends
to give a map $\alpha_1:Y_{\Sigma}/T(N_{m_0})\rightarrow {\bf C}\times
{\bf R}^{n-1}$. It is easy to check now that $\alpha_1$ is a homeomorphism
making the diagram (2.1) commute, and that $f=\im 
\circ \alpha_1\circ q_1$, as desired. $\bullet$
\bigskip

We next comment as to when the hypotheses of Theorem
2.2 can be achieved.

First, we recall the standard construction of $Y_{\Sigma}$ as a symplectic
quotient: see [3] for details. Here as always $Y_{\Sigma}$ is
assumed to be non-singular. Let $\Sigma(1)$ denote the set
of one-dimensional faces of $\sigma$. For each $\sigma\in\Sigma(1)$
there is a toric divisor $D_{\sigma}$ corresponding to $\sigma$, and
$\{D_{\sigma}|\sigma\in\Sigma(1)\}$ generates $\Pic Y_{\Sigma}$.
Define a map
$$\boldz^{\Sigma(1)}\mapright{\pi} N$$
taking a standard basis vector $e_\sigma$ of $\boldz^{\Sigma(1)}$,
$\sigma\in\Sigma(1)$, to the generator of the corresponding 
one-dimensional face of $\Sigma$. Let $K=\ker\pi$, and assume
(as will always be the case in our examples of interest) that
$\pi$ is surjective. For $I\subseteq\Sigma(1)$, define $e_I\subseteq
{\bf C}^{\Sigma(1)}$ by
$$\{(z_j)_{j\in\Sigma(1)}|\hbox{$z_i=0$ for $i\not\in I$}\}$$
and let $\bar I=\Sigma(1)\setminus I$. Let 
$$S=\{I\subseteq\Sigma(1)|\hbox
{$\bar I$ does not span a cone in $\Sigma$}\}.$$
Define $\U_{\Sigma}
\subseteq {\bf C}^{\Sigma(1)}$ by
$$\U_{\Sigma}={\bf C}^{\Sigma(1)}\setminus\bigcup_{I\in S} e_I.$$
$\U_{\Sigma}$ inherits the standard K\"ahler form 
$$\omega={i\over 2} \sum_{j\in\Sigma(1)} dz_j\wedge d\bar z_j.$$

Let ${\bf t}$ denote the Lie algebra of $T(\boldz^{\Sigma(1)})$
and ${\bf k}$ denote the Lie algebra of $T(K)$. There is a well-known
correspondence between elements of ${\bf t}^*$ and real divisors: an
element $\alpha\in {\bf t}^*$ corresponds to the real divisor
$$D_{\alpha}=\sum_{\sigma\in\Sigma(1)} \alpha(e_{\sigma}) D_{\sigma}.$$
On the other hand, such an $\alpha$ defines a piecewise linear function
on the support of the fan $\Sigma$. Indeed, for each cone $\tau$ of
$\Sigma$, choose $m_{\tau}\in M$ so that $\langle m_{\tau},\pi(e_{\sigma})
\rangle=-\alpha(e_{\sigma})$ for each one-dimensional face $\sigma$ of
$\tau$. Then define
$\phi_{\alpha}:|\Sigma|\rightarrow {\bf R}$ by $\phi_{\alpha}(x)=
\langle m_{\tau},x\rangle$ for $x\in\tau$. Then the ${\bf R}$-divisor
$D_{\alpha}$ is ample if and only if $\phi_{\alpha}$ is strictly upper convex,
i.e.
$$\langle m_{\tau},\pi(e_{\sigma})\rangle > -\alpha(e_\sigma)$$
whenever $\sigma$ is not a face of an $n$-dimensional cone $\tau$. 
Two divisors are linearly
equivalent if the corresponding piecewise linear functions differ by
a linear function; thus $\Pic Y_{\Sigma}\otimes {\bf R}\cong {\bf k}^*$
naturally via the projection $p:{\bf t}^*\rightarrow {\bf k}^*$.

Now the action of $T(\boldz^{\Sigma(1)})$ on $\U_{\Sigma}$ induces
the standard moment map $\nu:\U_{\Sigma}\rightarrow {\bf t}^*$ and
$\nu'=p\circ \nu:\U_{\Sigma}\rightarrow {\bf k}^*$ is the moment map
for the $T(K)$ action on $\U_{\Sigma}$.

Let $\alpha\in {\bf k}^*$ be the class of an ample ${\bf R}$-divisor
on $Y_{\Sigma}$. Then by Proposition 3.1.1 of [3], 
$\nu'^{-1}(\alpha)/T(K)$ is homeomorphic to $Y_{\Sigma}$, and the
induced, reduced symplectic form $\omega_{\alpha}$ on $Y_{\Sigma}$
has cohomology class equal to $D_{\alpha}$. This symplectic form
is a K\"ahler form.

The residual $T(N)$ action on $Y_{\Sigma}$ preserves this K\"ahler form,
so this gives an example of a $T(N)$-invariant K\"ahler form in each
K\"ahler class. Furthermore, this action induces the moment map
$\mu:Y_{\Sigma}\rightarrow p^{-1}(\alpha)$, the latter being a translation
of ${\bf n}^*$, where
${\bf n}$ is the Lie algebra of $T(N)$. Here we use the exact sequence
$$0\mapright{} {\bf n}^*\mapright{} {\bf t}^* \mapright{p} {\bf k}^*
\mapright{}0.$$
By [3], \S 3.2, the image of $\mu$ is $P_{\alpha}=
p^{-1}(\alpha)\cap \mu(\U_{\Sigma})$. It is then not difficult
to see that $P_{\alpha}$ is closed and
$\mu:Y_{\Sigma}\rightarrow {\bf n}^*$ is proper. Furthermore,
knowing the image of $\mu$ allows us to determine the discriminant
locus of the special Lagrangian fibration precisely.

{\it Example 2.3.} Let $N=\boldz^3+{1\over 3}(1,1,1)$, and let
$\sigma$ be the cone spanned by $(1,0,0)$, $(0,1,0)$ and $(0,0,1)$. $\Sigma$
is then obtained by subdividing $\sigma$ at ${1\over 3}(1,1,1)$, giving
$$\epsfbox{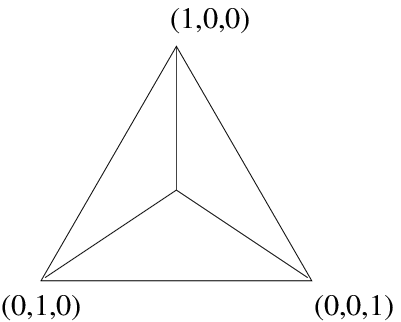}$$
Then 
$$\Sigma(1)=\{(1,0,0),(0,1,0),(0,0,1),{1\over 3}(1,1,1)\}$$
and the matrix for $\pi$ is
$$\pmatrix{1&0&0&1/3\cr
0&1&0&1/3\cr
0&0&1&1/3\cr}$$
with kernel generated by $(1,1,1,-3)$. Thus
$$\U_{\Sigma}={\bf C}^4\setminus \{z_1=z_2=z_3=0\}.$$
The moment map $\nu'$ is given by
$$\nu'(z_1,z_2,z_3,z_4)=|z_1|^2+|z_2|^2+|z_3|^2-3|z_4|^2.$$
For $a>0$, let $(Y_{\Sigma},\omega_{\Sigma})$ be given
by $\nu'^{-1}(a)/T(K)$. Now 
$$p^{-1}(a)=\{(r_1,r_2,r_3,(r_1+r_2+r_3-a)/3)|r_1,r_2,r_3\in {\bf R}\},$$
so $P_{a}$ can be identified with the set
$$\{(r_1,r_2,r_3)\in {\bf R}^3| r_1,r_2,r_3\ge 0, \quad r_1+r_2+r_3\ge a\}.$$
Finally, the discriminant locus of the induced special Lagrangian 
fibration $f=(g,\mu_0):Y_{\Sigma}\rightarrow{\bf R}^3$ 
is the planar graph which is the image of the 
1-skeleton of the boundary of $P_{a}$
under the projection ${\bf R}^3\rightarrow {\bf R}^2$
given by $(r_1,r_2,r_3)\mapsto (r_1-r_2,r_1-r_3)$. This is 
$$\epsfbox{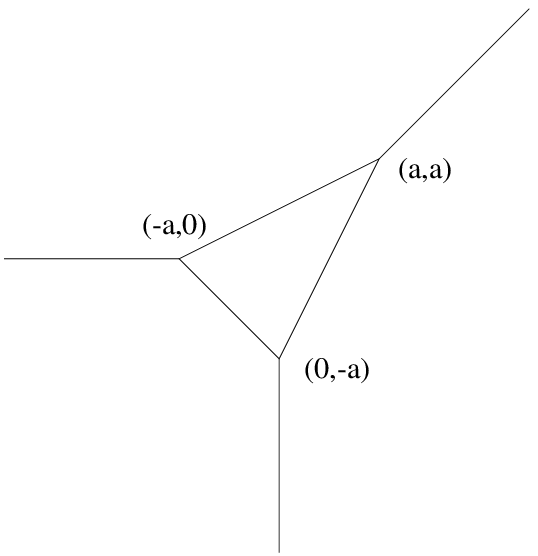}$$

\bigskip
In these examples, the moment map $\mu$ is always proper. More generally,
however, this need not be the case. Since the above metrics are not,
in general, Ricci-flat, we might be interested in a wider range of
metrics.
Now in any event
the action of $T(N)$ is induced by the action of $T_{\bf C}(N)$,
so it follows from [14], Convexity Theorem, \S 7, 
that the image of the moment map $\mu$ is
convex with connected fibres, so $\mu$ identifies $\mu(Y_{\Sigma})$
with $Y_{\Sigma}/T(N)$ in general. However, 
$\mu(Y_{\Sigma})$ may not be closed, in which case $\mu$ is not
proper. Thus we need some
additional asymptotic conditions on $\omega$.

For example in [34,35], existence of complete Ricci-flat metrics
on some non-compact manifolds was proven. More precise information in
certain cases was given in [19,20], where 
it was proved that there exists ALE or
quasiALE Ricci-flat metrics on 
crepant resolutions of ${\bf C}^n/G$, with $G\subseteq
SU(n)$. If $G$ is abelian, then ${\bf C}^n/G$ is toric. By the uniqueness
results of [19,20], these metrics are invariant under the
induced $T(N)$ action, and the above theorem applies. Furthermore,
the ALE or quasi-ALE conditions guarantee the moment map is proper, 
being asymptotic to the moment map with respect to the Euclidean metric.

In general, it is not known when Ricci-flat metrics exist on $Y_{\Sigma}$
for general Gorenstein cones $\sigma$. However, in some cases it is possible
to use the results of [34,35] to find further examples. In any event,
as long as the moment map with respect to $\omega$ is proper, its image
will coincide with the image of the moment map induced by $\omega'$
given by the symplectic reduction method above when the classes
$[\omega], [\omega']\in H^2(Y_{\Sigma},{\bf R})$ coincide. Thus we always
get a precise description of the discriminant locus. 
\medskip

Since we are interested in special Lagrangian fibrations because of the
SYZ conjecture, it is actually more interesting to construct proper
special Lagrangian fibrations. In [10], Remark 3.5, we noted we can
construct a topological ``properification'' of the map $f:Y_{\Sigma}
\rightarrow {\bf R}^3$. Here, we take an alternative route and construct
proper special Lagrangian fibrations as a generalization of [8],
Example 1.2.

\proclaim Theorem 2.4. Let $\omega$ be the K\"ahler form of a K\"ahler
metric on $Y_{\Sigma}$, invariant under the action of $T(N)$. Let $\Omega$
be the holomorphic $n$-form on $Y_{\Sigma}$ which restricts to 
$dz_1\wedge\cdots\wedge dz_n$ on $T_{\bf C}(N)$. Let $\mu:Y_{\Sigma}
\rightarrow {\bf R}^n$ be the moment map of the $T(N)$ action, and assume
$\mu$ is proper. Let $\mu_0:Y_{\Sigma}\rightarrow {\bf R}^{n-1}$ be the 
moment map associated to the $T(N_{m_0})$-action. Set
$$Y'_{\Sigma}=Y_{\Sigma}\setminus\{ 1+\prod_{i=1}^n z_i=0\}.$$
Here $\prod_{i=1}^n z_i$ is identified with the character $m_0$, and
as such, defines a regular function on $Y_{\Sigma}$, as we saw in the
proof of Theorem 2.2. Let
$$\Omega'={\Omega
\over i^n(1+\prod_{i=1}^n z_i)}$$
be a nowhere vanishing holomorphic $n$-form on $Y_{\Sigma}'$.
Then $f'=(\log |1+\prod_{i=1}^n z_i|,\mu_0):Y'_{\Sigma}
\rightarrow {\bf R}^n$ is a proper special Lagrangian fibration
with respect to $\omega$, $\Omega'$, with the same discriminant
locus as that of $f:Y_{\Sigma}\rightarrow {\bf R}^n$ constructed in
Theorem 2.2. Furthermore, the general fibre is an $n$-torus.

Proof. Let $\mu_0=(f_1,\ldots,f_{n-1})$ be as in the proof of Corollary
1.3. Then
$$\eqalign{\iota(X_{f_1},\ldots,X_{f_{n-1}})\Omega'
&=\pm{(-2i)^{n-1}d(z_1\cdot\cdots\cdots z_n)
\over i^n(1+\prod_{i=1}^n z_i)}\cr
&=\pm(-2)^{n-1} d(i\log(1+\prod_{i=1}^n z_i)),\cr}$$
so by Theorem 1.2, $$f'=(\im(i\log(1+\prod_{i=1}^n z_i)),\mu_0)
=(\log|1+\prod_{i=1}^n z_i|,\mu_0)$$
is a special Lagrangian fibration with respect to $\omega,\Omega'$.

Now consider $c\in{\bf R}^{n-1}$ and $\mu_0^{-1}(c)/T(N_{m_0})$,
with $\mu_0:Y_{\Sigma}\rightarrow {\bf R}^{n-1}$ (rather than its restriction
to $Y'_{\Sigma}$). Now in the proof of Theorem 2.2, it was
shown that $m_0$ (as a regular function) descended to 
$\mu_0^{-1}(c)/T(N_{m_0})$ to give an isomorphism $m_0:\mu_0^{-1}(c)/T(N_{m_0})
\rightarrow {\bf C}$. Then $m_0$ induces an isomorphism $m_0:(\mu_0^{-1}(c)
\cap Y'_{\Sigma})/T(N_{m_0})\rightarrow {\bf C}\setminus \{-1\}$.
Furthermore $\log |1+m_0|$ gives an $S^1$-fibration 
$(\mu_0^{-1}(c)\cap Y_{\Sigma}')/T(N_{m_0})
\rightarrow {\bf R}$. The inverse image of the general fibre of this
map in $\mu_0^{-1}(c)$ is then $T^{n-1}\times S^1=T^n$. Thus the general
fibre of $f'$ is an $n$-torus. Now orbits of $T(N_{m_0})$ only 
drop dimension when $m_0=0$. But when $m_0=0$, $\log |1+m_0|=0$,
and this makes it clear the discriminant locus coincides with that of
$f$. $\bullet$

{\it Remark 2.5.} (1) While $Y'_{\Sigma}$ is certainly not a 
holomorphic partial compactification of $Y_{\Sigma}$ (being contained
in $Y_{\Sigma}$), it is possible to prove $Y'_{\Sigma}\rightarrow
{\bf R}^n$
does coincide with a topological partial compactification of $Y_{\Sigma}
\rightarrow {\bf R}^n$.

(2) All singular fibres $f'^{-1}(b)$ of the above fibration have the same
basic structure: there is a fibration $f'^{-1}(b)
\rightarrow S_1$ with all but one 
fibre a $T^{n-1}$, with the remaining fibre
being a torus of dimension between $0$ and $n-2$.

(3) Unlike the non-proper case, we should not expect there to exist a complete
K\"ahler metric $\omega$ on $Y'_{\Sigma}$ satisfying
$\omega^n$ proportional to $\Omega'\wedge\bar\Omega'$.

\bigskip

{\hd \S 3. Deformations of toric singularities.}

The deformation theory of toric Gorenstein singularities is controlled
by the combinatorics of the corresponding cones. In particular, if
$Y_{\sigma}$ is an isolated toric singularity, then [1] gives a beautiful
description of the versal deformation space of $Y_{\sigma}$. Irreducible
components of this versal deformation space are in one-to-one 
correspondence with maximal Minkowski decompositions of
the polytope $P$ ($P$ as in \S 2, with $\sigma$ the cone over $P$).

Altmann's construction is as follows. Let $N$, $M$ and $\sigma$ be
as in \S 2, and assume $Y_{\sigma}$ has only an isolated singularity.
Now $\sigma$ is a cone over a polytope $P$ contained in the affine
hyperplane $\langle m_0,\cdot\rangle=1$. By choosing some element $n_0$
such that $\langle m_0,n_0\rangle=1$, we can identify  $P$
with $P-n_0$ in the hyperplane $L_{\bf R}\subseteq N_{\bf R}$
given by $L=m_0^{\perp}$. Let $P=R_0+\cdots+R_p$ be a Minkowski
decomposition of $P$ inside $L_{\bf R}$. What this means is that
$R_0,\ldots,R_p$ are convex subsets of $L$ such that 
$$P=\{r_0+\cdots+r_p|r_i\in R_i\}.$$

{\it Example 3.1.} We focus on the prettiest example, a cone over a del 
Pezzo surface of degree 6. We can take $\sigma$ to be generated by
$$n_1,\ldots,n_6=(0,0,1),(1,0,1),(2,1,1),(2,2,1),(1,2,1),(0,1,1)$$
so that $P$ in ${\bf R}^2$ is
$$\epsfbox{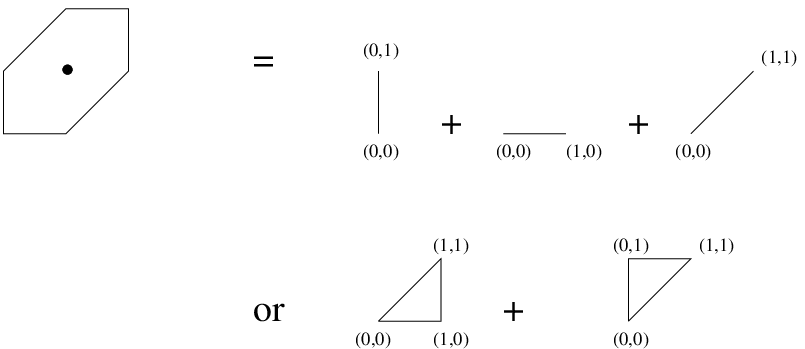}$$
with two different Minkowski decompositions.

\medskip
Now to each such Minkowski decomposition $P=R_0+\cdots+R_p$,
Altmann constructs a flat deformation of $Y_{\sigma}$ as follows. Let
$N'=L\oplus\boldz^{p+1}$, and let $e_0,\ldots,e_p$ denote the standard
basis of $\boldz^{p+1}$ and define
$$\tilde\sigma=Cone\left(\bigcup_{k=0}^p(R_k\times \{e_k\})\right)
\subseteq N'_{\bf R}$$
where $Cone(S)$ denotes the cone generated by the set $S\subseteq N'_{\bf R}$.
If one writes $N=L\oplus\boldz n_0$, there is a diagonal embedding
$N\hookrightarrow N'$ given by $l+an_0\mapsto l+a(e_0+\cdots+e_p)$.
Under this embedding $\sigma=\tilde\sigma\cap N_{\bf R}$, and hence we obtain
a closed embedding $Y_{\sigma}\hookrightarrow Y_{\tilde\sigma}$. On the other
hand, under the projection $N'\rightarrow \boldz^{p+1}$, $\tilde \sigma$
maps to the cone generated by $e_0,\ldots,e_p$, and this induces a morphism
$Y_{\tilde\sigma}\rightarrow {\bf C}^{p+1}$. Altmann proves the
composed morphism $f:Y_{\tilde\sigma}\rightarrow {\bf C}^{p+1}/{\bf C}
(1,\ldots,1)$ is a flat deformation of $Y_{\sigma}$, with $f^{-1}(0)=
Y_{\sigma}$. More explicitly, the surjection $N'\rightarrow\boldz^{p+1}$
gives an inclusion $\boldz^{p+1}\hookrightarrow M'$, with $e_0^*,
\ldots,e_p^*$ mapping to elements of $M'$ corresponding to 
characters $t_0,\ldots,t_p$. These characters extend to regular functions
on $Y_{\tilde\sigma}$, and $f$ is given by 
$(t_0-t_1,\ldots,t_0-t_p)$.

The main point for us then is that the functions $t_0,\ldots,t_p$ are
invariant under the action of $T_{\bf C}(L)\subseteq T_{\bf C}(N')$, and
thus $T_{\bf C}(L)$ acts on the fibres of $f$. This gives the desired
$T^{n-1}$-action on deformations of $Y_{\sigma}$.

What about a $T(L)$-invariant holomorphic $n$-form on the fibres of $f$?
Well note that $\tilde\sigma$ is a Gorenstein canonical cone;
if $m_0'=e_0^*+\cdots+e_p^*\in M'$, then all generators of $\tilde\sigma$
evaluate to 1 on $m_0'$. Thus there is a nowhere vanishing holomorphic
$n+p$-form $\Omega$ on the smooth part of $Y_{\tilde\sigma}$ (whose restriction
to $T_{\bf C}(N')\subseteq Y_{\tilde\sigma}$ is described by Proposition
2.1). We then have

\proclaim Proposition 3.2. If $t_0,\ldots,t_p$ are coordinates
on ${\bf C}^{p+1}$, and $\partial_{t_0},\ldots,\partial_{t_p}$ are lifts
of the corresponding vector fields to $Y_{\tilde\sigma}$, then for
$x\in {\bf C}^{p+1}/{\bf C}(1,\ldots,1)$,
$$\Omega_x=(\iota(\partial_{t_1},\ldots,\partial_{t_p})\Omega)|_{f^{-1}(x)}$$
is a well-defined nowhere vanishing holomorphic $n$-form on the non-singular
part of $f^{-1}(x)$, which we write as $Y^{ns}_{\tilde\sigma,x}$. 
In
addition, $\Omega_x$ is $T(L)$-invariant. Finally, if $z_1,\ldots,z_{n-1}$
are a basis of characters for $T_{\bf C}(L)$, then 
$z_1,\ldots,z_{n-1},t_0,\ldots,t_p$ form a basis of characters for
$T_{\bf C}(N')$ and 
$$\Omega|_{T_{\bf C}(N')}={dz_1\wedge\cdots\wedge dz_{n-1}\over
\prod z_i}\wedge dt_0\wedge\cdots\wedge dt_p$$
so up to sign
$$\Omega_x|_{T_{\bf C}(N)\cap f^{-1}(x)}=
{dz_1\wedge\cdots\wedge dz_{n-1}\over
\prod z_i}\wedge dt_0.$$

Proof. That $\Omega_x$ is well-defined, independent of the lifts of the 
$\partial_{t_i}$'s is standard, and since $\partial_{t_1},\ldots,\partial_{t_p}$
are linearly independent at a non-singular point of $f^{-1}(x)$,
$\Omega_x$ is non-zero.  Also, $\Omega$ is invariant under the action of
$T(m_0'^{\perp})$, and $L\subseteq m_0'^{\perp}$, so $\Omega$ is invariant
under $T(L)$. Since $t_1,\ldots,t_p$ are also invariant under $T(L)$,
so is $\Omega_x$. Finally, the explicit form
for $\Omega$ follows from Proposition 2.1 and the explicit value for
$m_0'$. $\bullet$

\proclaim Proposition 3.3. If $\omega$ is a $T(L)$-invariant K\"ahler
form on $Y^{ns}_{\tilde\sigma,x}$, let $\mu:Y^{ns}_{\tilde\sigma,x}
\rightarrow {L_{\bf R}^*}$ be the moment map associated to the $T(L)$-action.
Then $f:Y^{ns}_{\tilde\sigma,x}\rightarrow{\bf R}\times L_{\bf R}^*$
given by $f=(\im(i^{n+1}t_0),\mu)$ is a special Lagrangian fibration.
Furthermore, $f$ is surjective if $\mu$ is and the general fibre is
diffeomorphic to ${\bf R}\times T(L)$. If $x$ is represented by
$(x_0,\ldots,x_p)\in {\bf C}^{p+1}$, then the discriminant locus is contained
in the union of $p+1$ hyperplanes $\{\im(i^{n+1}(x_0-x_k))| k=0,\ldots,p\}
\times L_{\bf R}^*$. 

Proof. The fact that $f$ is special Lagrangian follows immediately
from the form of $\Omega_x$ given in Proposition 3.1 and the same type
of calculation as performed in Corollary 1.3. Now if $y\in Y^{ns}_{\tilde
\sigma,x}$ then $\mu^{-1}(\mu(y))/T(L)$ is isomorphic to the 
categorical quotient $Y^{ns}_{\tilde\sigma,x}//T_{\bf C}(L)$, which is
isomorphic to ${\bf C}$ with holomorphic coordinate $t_0$.
Thus the fibre $f^{-1}(f(y))$ is an inverse image of a straight line
in ${\bf C}$ under the quotient map $\mu^{-1}(\mu(y))\rightarrow
{\bf C}$. Thus the general fibre is $T(L)\times {\bf R}$.
Clearly also $f$ is surjective if $\mu$ is. 

The discriminant locus $\Delta$ is the image of the union of $T(L)$
orbits of dimension $<n-1$. Now a subcone $\tau$ of dimension $k$
of $\tilde\sigma$ corresponds to a codimension $k$ orbit of $T_{\bf C}(N')$
which is fixed by $T_{\bf C}({\bf R}\tau\cap N')$. Thus the $T(L)$ orbits
on this stratum drop dimension if ${\bf R}\tau\cap L_{\bf R}\not=0$.
The one-dimensional faces of $\tilde\sigma$ are generated by
$n\times e_k$ where $n$ is a vertex of the polytope $R_k$. Thus
a face $\tau$ of $\tilde\sigma$ has ${\bf R}\tau\cap L_{\bf R}\not=0$
if and only if it contains two one-dimensional faces generated by
$n_1\times e_k$ and $n_2\times e_k$ for some $k$, for $n_1,n_2$ two 
vertices of $R_k$. Necessarily $n_1$ and $n_2$ are the endpoints of
an edge of $R_k$. Thus all minimal faces $\tau$ such that ${\bf R}\tau
\cap L_{\bf R}\not= 0$ are two-dimensional faces spanned by $n_1\times e_k$
and $n_2\times e_k$. Now the function $t_k$ is necessarily zero on
the corresponding codimension 2 stratum. If this codimension 2 stratum
is called $D_{\tau}\subseteq Y_{\tilde\sigma}$, then on $D_{\tau}\cap
Y_{\tilde\sigma,x}$, $t_k-t_0=x_k-x_0$ so $t_0=x_0-x_k$. Thus
$f(D_{\tau}\cap Y_{\tilde\sigma,x})$ is contained in the hyperplane
given by $\{\im(i^{n+1}(x_0-x_k))\}\times L_{\bf R}^*$. $\bullet$

{\it Example 3.4.} Continuing with Example 3.1, the discriminant
locus depends on the choice of the two decompositions. For general choice
of $x$, the discriminant locus in the first splitting is contained
in 3 different planes. There are three choices of two-dimensional $\tau$
yielding components of the discriminant locus, and for each $\tau$,
$D_{\tau}\cap Y_{\tilde\sigma,x}$ (for general $x$, $Y_{\tilde\sigma,x}$
is already non-singular) consists just of a ${\bf C}^*$. The image of
this ${\bf C}^*$ under the moment map is a straight line, and the fibres
of this map are connected (by [14]). Depending on the properties
of $\mu$, this image is either a line segment, a ray, or a line infinite
in both directions. So $\Delta$ looks like
$$\epsfbox{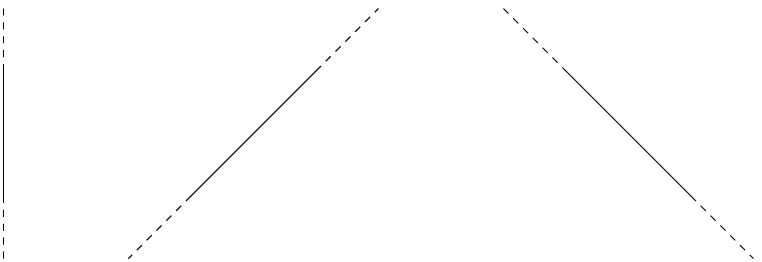}$$
where each line is in a parallel plane. As $x\rightarrow 0$, these planes
will converge to the same plane, producing, for suitable choice of 
$\omega$, a discriminant locus for $x=0$ of
$$\epsfbox{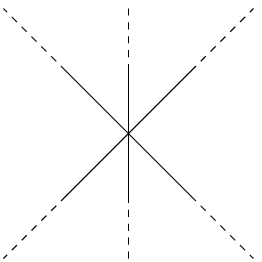}$$
If instead we choose the second smoothing, then similar arguments show
that $\Delta$ looks like
$$\epsfbox{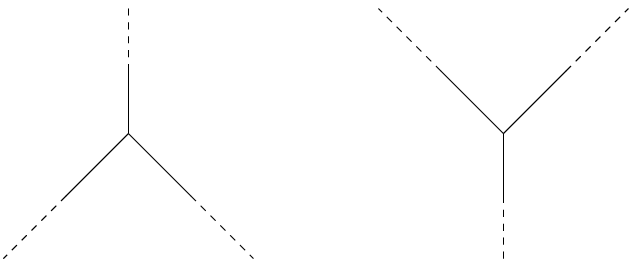}$$
Again, we get the same picture as above as $x\rightarrow 0$. 

If we take a crepant resolution of $Y_{\tilde\sigma,0}=Y_{\sigma}$, then
by \S 2 we obtain a discriminant locus which looks like
$$\epsfbox{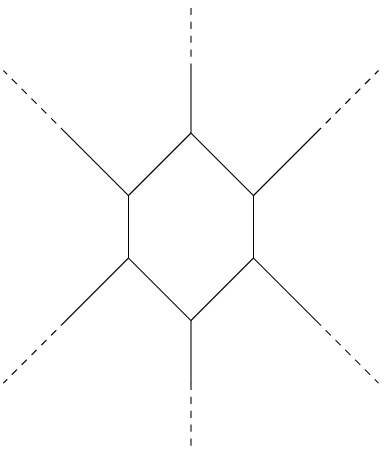}$$

Another, simpler, example, is the ordinary double point, given by
$$\epsfbox{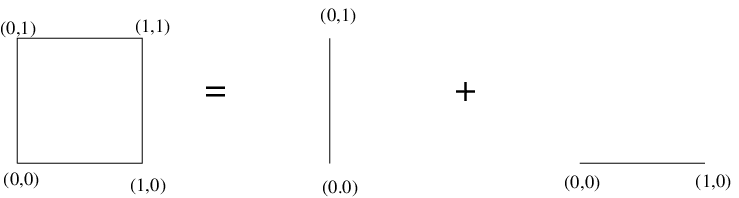}$$
The discriminant locus of the smoothing consists of two lines in
different planes. As $x\rightarrow 0$ these planes converge, and then
there are two different small resolutions of the ordinary double point. This
gives a picture
$$\epsfbox{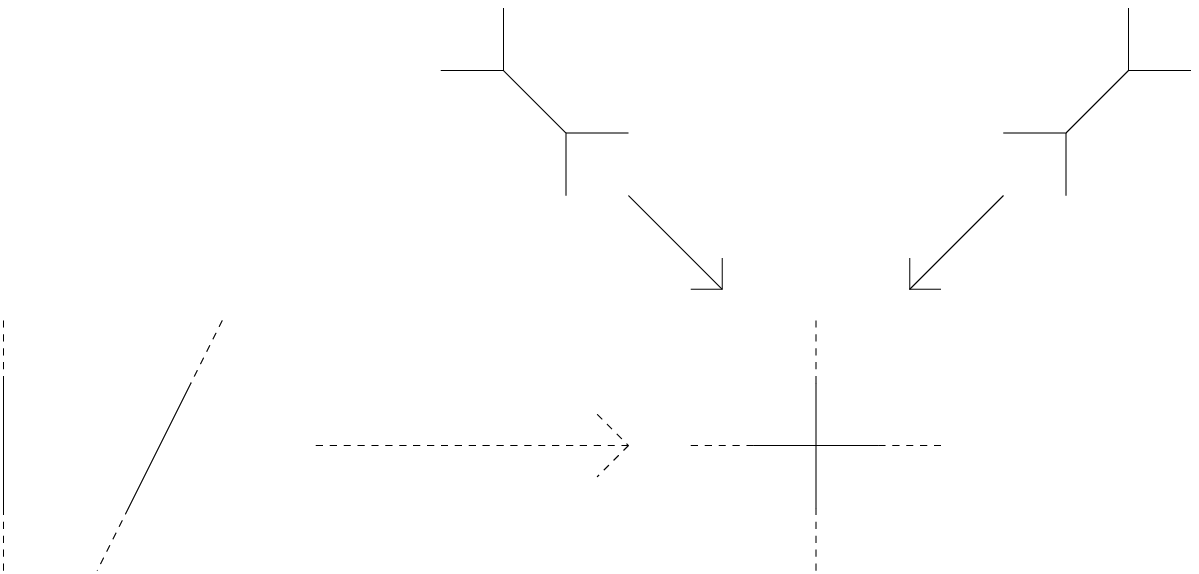}$$
which shows how the discriminant locus changes under smoothing and resolution.

We have not discussed the choice of the metrics on these smoothings,
but in some examples one can find Ricci-flat metrics. In the ordinary
double point case, Stenzel [32] has constructed an explicit Ricci-flat
metric, while in some examples [34,35] apply.

\bigskip
{\hd \S 4. Local mirror symmetry and 
connections with the work of Ruan and Joyce.}

Let us begin this more informal discussion by asking the question: 
how do we construct the mirror to $Y_{\Sigma}$, where $Y_{\Sigma}$ is
as in \S 2? We shall focus on the three-dimensional case. This should make
sense in the context of local mirror symmetry as developed in [4]. (See
also related forms of local mirror symmetry in [5], [17], and [18]).
If we follow the SYZ philosophy, then we would need to construct duals
of the special Lagrangian fibrations $f:Y_{\Sigma}\rightarrow {\bf R}^3$
of Theorem 2.2. The difficulty is that the general fibre of $f$
is $T^2\times {\bf R}$, which we can't dualize. However, we have at
least constructed, via Theorem 2.4, a topological
``properification'' $f':Y'_{\Sigma}\rightarrow {\bf R}^3$
of $f$. (We will ignore the metric properties of $Y_{\Sigma}$ and 
$Y'_{\Sigma}$ for the moment). We can then dualize $f'$ topologically.
Indeed, from the description of the singular fibres in Remark 2.5, (2),
as well as the calculation of monodromy in [8], Example 1.2, we can
see that
$f'$ has only semi-stable fibres, and in fact satisfies the hypotheses
of Corollary 2.2 of [10]. Thus a topological dual $\check f:\check Y_{\Sigma}
\rightarrow {\bf R}^3$ of $f':Y'_{\Sigma}\rightarrow
{\bf R}^3$ exists.

However, it is worthwhile describing this dual explicitly. 
First, observe that the discriminant locus $\Delta$
of $f'$ is a trivalent graph (homeomorphic to the 1-skeleton
of $\partial\mu(Y_{\Sigma})$). The fibres over the edges of the
graph are of type $(2,2)$ in the notation of [9] and [10] (i.e.
a product of a circle with a Kodaira type $I_1$ fibre) and type
$(1,2)$ at all vertices\footnote{*}{For generic good $T^3$-fibrations,
which only have singular fibres of type $(2,2)$, $(2,1)$ and $(1,2)$,
I am going to second Dave Morrison's suggestion that they be called
generic, negative and positive singular fibres respectively, the words
negative and positive referring to the sign of the Euler characteristic
of the fibre. W.-D. Ruan introduced his own notation of type II and III
for these singular fibres.}.
Thus the dual $\check f:\check Y_{\Sigma}\rightarrow {\bf R}^3$ has only
type $(2,2)$ and type $(2,1)$ fibres.

Let's describe the monodromy of the fibration $f'$.
Fix a basis $e_1,e_2,e_3$
of $N$ so that that $m_0=e_1^*+e_2^*+e_3^*$ and take $f_1=e_1-e_2, f_2=e_1
-e_3$ to be a basis for $N_{m_0}$. If we take any point $b\in {\bf R}^3
\setminus\Delta$, then the fibre $f'^{-1}(b)$ has a $T(N_{m_0})$-action
which allows us to identify $N_{m_0}$ with a sublattice of
$H_1(f'^{-1}(b),\boldz)$. We can then choose an element $f_3\in 
H_1(f'^{-1}(b),\boldz)$ such that $f_1,f_2,f_3$ form a basis for
$H_1(f'^{-1}(b),\boldz)$. It is clear $f_1$ and $f_2$ will be monodromy
invariant 1-cycles.

An edge $l$ of $\Delta$ is the image under $f'$ of a codimension-two
$T_{\bf C}(N)$-orbit, which in turn corresponds to a dimension 2 face
$\tau$ of $\Sigma$ with generators $n_1,n_2$.
It then follows from [10], Proposition 3.3 and
Example 2.8, that if $b$ is chosen near the edge $l$, $\gamma:
S^1\rightarrow {\bf R}^3\setminus\Delta$ a suitably oriented simple
loop about $l$ based at $b$, and $n_1-n_2=a_1f_1+a_2f_2$, then the
monodromy transformation $T:H_1(f'^{-1}(b),\boldz)\rightarrow H_1(f'^{-1}(b),
\boldz)$ about $\gamma$ is, in the basis $f_1,f_2,f_3$, 
$$T=\pmatrix{1&0&a_1\cr 0&1&a_2\cr 0&0&1\cr}.$$
The topological dual $\check f:\check Y_{\Sigma}\rightarrow {\bf R}^3$
is constructed as follows. Let $N_{m_0}^*=\hom_{\boldz}(N_{m_0},\boldz)$.
Set $\bar X=T(N_{m_0}^*)\times {\bf R}^3$,
and let $S\subseteq X$ be a topological surface constructed as follows.
$S$ projects to $\Delta\subseteq {\bf R}^3$, and  for each edge
$l$ of $\Delta$, $(T(N_{m_0}^*)\times l)
\cap S$ is a cylinder fibering in circles over $l$, with the
circle homotopic to $T((n_1-n_2)^{\perp})\subseteq T(N_{m_0}^*)$. At a
vertex, which corresponds to a two-dimensional face of $\Sigma$,
spanned by $n_1,n_2$ and $n_3$, the classes $n_1-n_2,
n_2-n_3$ and $n_3-n_1$ add to zero, and thus the three cylinders
can be glued above the vertex. This gives a topological manifold $S$.

Now let $X=\bar X\setminus S$, $\pi:Y\rightarrow X$ a principal
$S^1$-bundle with
Chern class $(0,\pm 1)\in H^2(X,\boldz)=H^2(\bar X,\boldz)\oplus\boldz$. 
Then there is a topological manifold $\bar Y$ containing $Y$ and a diagram
$$\matrix{Y&\hookrightarrow&\bar Y\cr
\mapdown{\pi}&&\mapdown{\bar\pi}\cr
X&\hookrightarrow&\bar X\cr}$$
such that $\bar\pi$ is proper and the $S^1$-action on $Y$ extends
to an $S^1$-action on $\bar Y$, with $\bar\pi^{-1}(S)\cong S$
(see [10], Proposition 2.5). Taking $\check Y_{\Sigma}=\bar Y$ and
$\check f$ the composition $\bar Y\rightarrow\bar X\rightarrow {\bf R}^3$,
we obtain the topological dual of $f':Y'_{\Sigma}\rightarrow {\bf R}^3$.
These are dual in the sense that the monodromy representations are dual,
which is the only topological measure of duality.

If we are only interested in the topology, this is the end of the
story. But to get further insight into the picture, let us consider
the local mirror symmetry picture of [4].

One way to interpret the suggestions of [4] is as follows: the mirror
to $Y_{\Sigma}$ (not $Y_{\Sigma}'$!) is a curve $C\subseteq 
({\bf C}^*)^2=T_{\bf C}(N_{m_0}^*)$ 
whose Newton polygon is the translation of the
polygon $P$ into the plane $N_{m_0}\otimes {\bf R}$.
Calculations of certain period integrals
should yield the mirror map and predictions for Gromov-Witten invariants
on $Y_{\Sigma}$. These period integrals are integrals of ${dz_1\wedge dz_2
\over z_1z_2}$ on $({\bf C}^*)^2$ over 2-cycles with boundary on $C$.
Such integrals satisfy standard Picard-Fuchs equations.

The basic claim is that the pair $S\subseteq T(N_{m_0}^*)\times {\bf R}^2$
contained in $\bar X$ (where ${\bf R}^2\subseteq {\bf R}^3$ is the plane 
containing $\Delta$) is the same as the pair $C\subseteq T_{\bf C}(N_{m_0}^*)$.
This will make the connection both between this circle of ideas and
local mirror symmetry, as well as the connection with Ruan's work, clear.

To make this connection, we need to be more precise about our choice of
the equation for $C$. First, choose a K\"ahler class on $Y_{\Sigma}$.
(The choice of K\"ahler class determines $\Delta\subseteq {\bf R}^3$.)
As in \S 2, this can be thought of as a strictly convex function
$\phi:|\Sigma|\rightarrow {\bf R}$, which we can restrict to the polygon
$P$. With coordinates $z_1,z_2$ on $({\bf C}^*)^2$ corresponding to
the basis $f_1,f_2$ of $N_{m_0}$, we can consider the family of curves
$C_t$ given by the equation $h_t=0$, $t>0$, where
$$h_t=\sum_{(a,b)\in P\cap N_{m_0}} t^{\phi(a,b)}m_{a,b} z_1^a z_2^b,$$
with $m_{a,b}\in {\bf C}$.
Let $\nu:T_{\bf C}(N_{m_0}^*)\rightarrow {\bf R}^2$ be the moment map 
$\nu(z_1,z_2)=(\log |z_1|,\log |z_2|)$. 

The following theorem is implicit in Ruan's work ([29,30])
and can also be proved using ideas of Viro [36] and Mikhalkin [27].

\proclaim Theorem 4.1. For $|t|$ close to zero, $\nu(C_t)$ is a fattening
of the graph $\Delta$, and there is a $C^0$-isotopy of $T_{\bf C}(N_{m_0}^*)
=T(N_{m_0}^*)\times {\bf R}^2$ identifying $C_t$ and $S$.

Ruan uses this to construct torus fibrations on toric hypersurfaces.
Thus he was led to his pictures of the discriminant loci of Lagrangian
fibrations by looking at the {\it dual} picture to the fibrations
$Y_{\Sigma}\rightarrow {\bf R}^3$ developed in [10] and here. These two
points of view complement each other nicely.

{\it Example 4.2.} Let us continue with Example 2.3. Here $m_0=(1,1,1)$,
and we can take a basis $f_1,f_2$ of $N_{m_0}$ with $f_1=(2/3,-1/3,-1/3)$,
$f_2=(-1/3,2/3,-1/3)$, so that if $P$ is translated to $N_{m_0}$, we
can take it to be the convex hull of $f_1,f_2$ and $-f_1-f_2$. The set
of integral points of $P$ then corresponds to the monomials $z_1,z_2,1$
and $z_1^{-1}z_2^{-1}$. We can take
$$h_t=t(z_1+z_2+z_1^{-1}z_2^{-1})+1,$$ in which case, for $t$ small,
the image of $C_t$ under $\nu$ looks like
$$\epsfbox{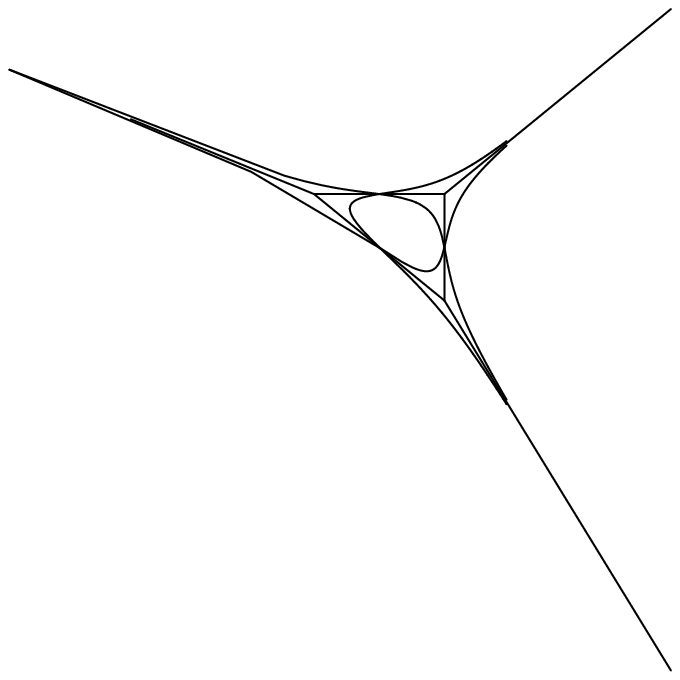}$$
We have in fact only drawn the boundary of $\nu(C_t)$, and superimposed
the codimension two discriminant locus $\Delta$ of $f:Y_{\Sigma}
\rightarrow {\bf R}^3$ with symplectic form given by the monomial-divisor
mirror map. This discriminant locus $\Delta$ is lying in the interior
of $\nu(C_t)$. (The actual shape of $\Delta$ disagrees with the one drawn
in Example 2.3, because we have used a different basis for $N_{m_0}$.)

\bigskip
Now, so far we have produced a topological fibration $\check f:
\check Y_{\Sigma}\rightarrow {\bf R}^3$. How might we construct a special
Lagrangian fibration? To do this we must first realise $\check Y_{\Sigma}$
as an (almost) Calabi-Yau manifold. We should expect this structure to be
invariant under an $S^1$ action with fixed locus isomorphic to $C_t$.
Furthermore, we might expect this $S^1$ action to extend to a
${\bf C}^*$ action on $\check Y_{\Sigma}$. One way to accomplish this is as
follows. Let $h$ be a regular function on $({\bf C}^*)^2$, and let
$$Y_h=\{(x,y,z_1,z_2)\in {\bf C}^2\times ({\bf C}^*)^2|
xy=h(z_1,z_2)\}.$$
$Y_h$ has a ${\bf C}^*$-action given by $(x,y,z_1,z_2)\mapsto
(\lambda x,\lambda^{-1}y,z_1,z_2)$ for $\lambda\in {\bf C}^*$,
and the fixed locus is the curve $x=y=h=0$. We also need to choose a
K\"ahler and holomorphic 3-form on $Y_h$. The holomorphic form will
be
$$\Omega=i{dx\wedge dz_1\wedge dz_2\over xz_1z_2}
=-i{dy\wedge dz_1\wedge dz_2\over yz_1z_2}$$
on $Y_h$. We have more choice for $\omega$, but for convenience we will
take
$$\omega={i\over 2}\left(dx\wedge d\bar x+dy\wedge d\bar y
+{dz_1\wedge d\bar z_1\over |z_1|^2}+{dz_2\wedge d\bar z_2\over
|z_2|^2}\right).$$
(The form of the part of $\omega$ which is a pull-back from $({\bf C}^*)^2$
is crucial.)

Now try to construct a special Lagrangian fibration on $Y_h$ using
Theorem 1.2. The moment map of the $S^1$ action is
$\mu:Y_h\rightarrow
{\bf R}$ given by $\mu(x,y,z_1,z_2)=|x|^2-|y|^2$. 
We need, for each $c\in {\bf R}$, to find a special Lagrangian
fibration on $\mu^{-1}(c)/S^1$. Now $\mu^{-1}(c)/S^1$ is canonically isomorphic
to $({\bf C}^*)^2$ as a complex manifold, and 
$$\eqalign{\Omega_{red}&=\iota(2i(y\partial_y-x\partial_x))\Omega\cr
&=2{dz_1\wedge dz_2\over z_1z_2},\cr}$$
while $\omega_{red}$ can be calculated with some effort to be
$$\omega_{red}={i\over 2}\left( {dz_1\wedge d\bar z_1\over |z_1|^2}
+{dz_2\wedge d\bar z_2\over |z_2|^2} +{1\over\sqrt{c^2+4|h|^2}} dh\wedge
d\bar h\right).$$
We need to find a special Lagrangian fibration on $(({\bf C}^*)^2,
\omega_{red},\Omega_{red})$ depending on the value of $c$. Now if
we took the limit $c\rightarrow\infty$, we get to the case where 
$$\omega_{red}={i\over 2}\left( {dz_1\wedge d\bar z_1\over |z_1|^2}
+{dz_2\wedge d\bar z_2\over |z_2|^2}\right),$$
and the moment map $\nu:({\bf C}^*)^2\rightarrow {\bf R}^2$ is then a special
Lagrangian fibration. Thus the special Lagrangian fibrations for finite
$c$, if they exist, should be viewed as a deformation of $\nu$. For finite
$c\not=0$ one might hope to prove the existence of such a deformation using
pseudo-holomorphic curve techniques and Gromov compactness. However,
at $c=0$ the almost Calabi-Yau structure becomes singular, and it is 
difficult to predict the structure of a special Lagrangian fibration on
$Y_h$, if it should exist. The conjectural picture however is that there
exists a special Lagrangian fibration $f:Y_h\rightarrow {\bf R}^3$ given by
$(|x|^2-|y|^2,g)$, where $g:Y_h\rightarrow {\bf R}^2$ induces a special
Lagrangian fibration on each reduced space. Furthermore the
discriminant locus of $f$ will be contained in the hyperplane in
${\bf R}^3$ where the first coordinate is zero, and would be the image
of $x=y=h=0$ under the map $f$, i.e. the image of $h=0$ in
$({\bf C}^*)^2$ under some deformation of the moment map $\nu$.

Of course, if we are interested in a dual to $f':Y_{\Sigma}'
\rightarrow {\bf R}^3$ which is special Lagrangian, we should
take $h=h_t$, identifying $\check Y_{\Sigma}$ with $Y_{h_t}$ via
Thereom 4.1. Getting more speculative, we recall that dualizing should
exchange information at a deeper level, i.e. interchange the data
of the symplectic structure on $Y_{\Sigma}$ with the complex structure
on $\check Y_{\Sigma}$. The precise correspondence is understood
at the level of the mirror map, as calculated in [4]. I don't want
to go into details here, but the main point is that one wants to
compute periods of the holomorphic 3-form $\Omega$ on $\check Y_{\Sigma}$.
Taking $\check Y_{\Sigma}=Y_{h_t}$, note that $H_3(\check Y_{\Sigma},\boldz)
=H_3(\mu^{-1}(0),\boldz)$ (as gradient flow gives a retraction of
$\check Y_{\Sigma}$ onto $\mu^{-1}(0)$), and in turn $H_3(\mu^{-1}(0),\boldz)
\cong H_2(({\bf C}^*)^2,C_t,\boldz)$. In other words a 3-cycle
is the inverse image under the quotient map of a 2-chain with boundary
in $C_t$. Integrating over such a 3-cycle reduces to integrating
${dz_1\wedge dz_2\over z_1z_2}$ over such a 2-chain. These period integrals
satisfy the relevant Picard-Fuchs equations and define the right
mirror map as explored in [4]. This perhaps explains why $C_t
\subseteq ({\bf C}^*)^2$ should be considered as the mirror of $Y_{\Sigma}$.

As one further intuitive observation in this direction, let us try to 
explain why it makes sense to consider $C_t\subseteq ({\bf C}^*)^2$ to be
the mirror of $Y_{\Sigma}$ rather than $Y_{\Sigma}'$. $T$-duality,
to first approximation, should exchange long and short distances in the
fibres. Thus the fibres of $f:Y_{\Sigma}\rightarrow {\bf R}^3$ should be
viewed as a limit of tori with greater and greater radius in one direction,
so in the limit the fibre is ${\bf R}\times T^2$ rather than $T^3$.
The $T$-dual fibres should have this radius approaching zero. Intuitively,
it then appears natural to divide $Y_{h_t}$ by the $S^1$-action, as
the $S^1$'s should correspond to the ``small'' direction.
\medskip

This discussion should not be taken too seriously. It is clear that
there is much to understand, but I believe this circle of ideas and
examples will prove to be an excellent laboratory for exploring 
the more intricate questions surrounding the SYZ conjecture.

Perhaps the most pressing question is the following.
It was originally my and many others' hope that special Lagrangian
fibrations would be reasonably differentiable, and differentiability
implies certain conditions on the discriminant locus. For example,
the discriminant locus is Hausdorff codimension 2 if the map is 
$C^{\infty}$ (see [10], \S 1). 
However, Joyce has now given in [22] examples
of special Lagrangian fibrations which are only piecewise differentiable,
and whose discriminant locus is codimension one. In fact, his basic 
example is of a very similar flavour to the $S^1$-invariant setup above.
From these examples, I believe the likelihood is that if $f:Y_h\rightarrow
{\bf R}^3$ exists, it is only piecewise smooth, and the discriminant
locus is amoeba-like rather than a graph. 
Thus if there is a special Lagrangian fibration
$f:Y_{h_t}\rightarrow {\bf R}^3$, it is only a perturbation of the
topological fibration $\check f:\check Y_{\Sigma}\rightarrow {\bf R}^3$.
\bigskip

{\hd \S 5. The Future of the SYZ Conjecture.}

It is clear that Joyce's picture forces us to reconsider the full
strength version of the SYZ conjecture, as opposed to the topological
ones considered in [10] and [29]. The notion of dualizing topological
torus fibrations as developed in [10] will not be the right one. We would
expect that if a mirror pair $X,\check X$ possess special Lagrangian
fibrations $f:X\rightarrow B$, $\check f:\check X\rightarrow B$,
they will have different amoeba-like discriminant loci $\Delta$ and
$\check \Delta$. They will presumably be dual only in a relatively
crude topological sense in that the monodromy representations
$\rho:\pi_1(B\setminus\Delta)\rightarrow SL_3(\boldz)$ and 
$\check\rho:\pi_1(B\setminus\check\Delta)\rightarrow SL_3(\boldz)$ are dual
representations.

It is my current belief that the SYZ conjecture will make most sense in
a limiting picture. First let us recall from [16] certain structures
which appear naturally on moduli spaces of special Lagrangian submanifolds.
Let $(X,\omega,\Omega)$ be a Calabi-Yau manifold, $B$ a moduli space
of deformations of some special Lagrangian submanifold on $X$, along
with a universal family
$$\matrix{\U&\hookrightarrow&X\times B\cr
\mapdown{f}&&\cr
B&&\cr}$$
Let $p:\U\rightarrow X$ be the projection. All fibres of $f$ are
assumed to be smooth submanifolds of $X$. Then from [26],
we know $B$ is smooth, with a canonical identification of $\T_{B,b}$
with $\H^1(f^{-1}(b),{\bf R})$, the space of ${\bf R}$-valued harmonic 
one-forms,
on $f^{-1}(b)$. This is via the map $v\in\T_{B,b}\mapsto \iota(v)p^*\omega$,
where $v$ is pulled back to a vector field normal to $f^{-1}(b)$ in $\U$.

There are two important structures on $B$:
\item{(1)} an integral affine structure. If $U\subseteq B$ is a 
contractible open set, with coordinates $t_1,\ldots,t_n$, let 
$\gamma_1,\ldots,\gamma_n\in H_1(f^{-1}(b),\boldz)$ be a basis for
first homology varying continuously with $b$. Then the $1$-forms $\alpha_i$
given by
$$\partial/\partial t_j\mapsto \int_{\gamma_i} \iota(\partial/\partial t_j)
p^*\omega$$
on $B$ are closed and linearly independent ([16], Proposition 1). Thus there
exists a coordinate system $y_1,\ldots,y_n$ on $U$ with $\alpha_i=dy_i$,
and these coordinates are well-defined up to integral affine
transformations (elements of ${\bf R}^n\rtimes GL_n(\boldz)$). This
defines an integral affine structure on $B$.
\item{(2)} There is a metric on $B$ (the McLean metric) given by
$$g(\partial/\partial t_i,\partial/\partial t_j)=
-\int_{f^{-1}(b)}\iota(\partial/\partial t_i)p^*\omega \wedge
\iota(\partial/\partial t_j)p^*\im\Omega.$$
(One might want to normalize this metric in various ways).

Hitchin showed there is a compatability between the metric and affine
structure: locally there exists a function $K$ such that $g(\partial/\partial
y_i,\partial/\partial y_j)=\partial^2 K/\partial y_i\partial y_j$.
Kontsevich and Soibelman [23]
call this structure on $B_0$ of affine structure plus metric of this form
an {\it affine K\"ahler} (AK) manifold. Such a manifold was called {\it
Hessian} in earlier work of H. Shima: see [31] and references therein.
If in addition the function $K$ satisfies the real Monge-Amp\`ere equation
$\det \partial^2K/\partial y_i\partial y_j=constant$, Kontsevich and
Soibelman call such an Hessian manifold a {\it Monge-Amp\`ere manifold}.

Next recall the definition of Gromov-Hausdorff convergence.

\proclaim Definition 5.1. Let $(X,d_X)$, $(Y,d_Y)$ be two compact
metric spaces. Suppose there exists maps $f:X\rightarrow Y$
and $g:Y\rightarrow X$ (not necessarily continuous) such that
for all $x_1,x_2\in X$,
$$|d_X(x_1,x_2)-d_Y(f(x_1),f(x_2))|<\epsilon$$
and for all $x\in X$,
$$d_X(x,g\circ f(x))<\epsilon,$$
and the two symmetric properties for $Y$ hold. Then we say the
Gromov--Hausdorff distance between $X$ and $Y$ is at most $\epsilon$.
The Gromov--Hausdorff distance $d_{GH}(X,Y)$ is the infinum of all
such $\epsilon$.

There are two distinct situations we might want to apply this notion.
The first was discussed independently by myself and Wilson [11]
and Kontsevich and Soibelman in [23]. The second situation follows
naturally from these ideas. First, let
$\X\rightarrow\Delta$ be a flat family of degenerating Calabi-Yau
$n$-folds, with $0\in\Delta$ a large complex structure limit point
(or maximally unipotent boundary point). Let $t_i\in\Delta$ be a sequence
of points converging to $0\in\Delta$, 
and let $g_i$ on $\X_{t_i}$ be a Ricci-flat
metric normalized so that $Diam(\X_{t_i},g_i)$ remains constant.
Then general results about Gromov-Hausdorff convergence tell
us that a subsequence of $(\X_{t_i},g_i)$ converges to a metric space
$(X_{\infty},g_{\infty})$. 

Second, consider another sequence of metric spaces, whose existence (or
rather non-emptiness) is currently conjectural. Suppose that for
$t_i$ sufficiently close to $0$, there is a special Lagrangian $T^n$
whose homology class is invariant under monodromy near $0$. (This
is a property we expect to find of fibres of a special Lagrangian
fibration associated to a large complex structure limit point).
Let $B_{0,i}$ be the moduli space of deformations of this torus,
every point of $B_{0,i}$ corresponding to a smooth torus in
$\X_{t_i}$. The manifold $B_{0,i}$ comes equipped with the McLean
metric. We should then compactify $B_{0,i}\subseteq B_i$ in some
manner: probably taking the closure of $B_{0,i}$ in the space of
special Lagrangian currents on $\X_{t_i}$ is the right thing to do.
This should give a series of metric spaces $(B_i,d_i)$, which again,
if the McLean metric is normalized properly to keep the diameter
constant, may have a convergent subsequence,
converging to a compact metric space $(B_{\infty},d_{\infty})$.

The following is a slight souping up of the conjectures in [11]
and [23].

\proclaim Conjecture 5.2. If $(\X_{t_i},g_i)$ converges to $(X_{\infty},
g_{\infty})$ and $(B_i,d_i)$ is non-empty for large $i$ and converges
to $(B_{\infty},d_{\infty})$, then $B_{\infty}$ and $X_{\infty}$ are
isometric up to scaling. Furthermore, there is a subspace $B_0\subseteq 
B_{\infty}$ with $\Delta=B_{\infty}\setminus B_0$ of Hausdorff codimension
2 in $B_{\infty}$ such that $B_0$ is a
Monge-Amp\`ere manifold, with the
metric inducing $d_{\infty}$ on $B_0$.

This is a considerably weaker conjecture than the original full-strength
SYZ proposal on the existence of special Lagrangian fibrations. But following
the philosophy of Kontsevich and Soibelman, this should be sufficient
for most purposes.

{\it Remarks 5.2.} 
(1) This conjecture doesn't assume the existence of special Lagrangian
fibrations on $\X_{t_i}$ for any $i$. It would of course be nice if
this is the case, but taking Joyce's philosophy seriously means we only
see the codimension 2 structure in the limit. We expect that as $i\rightarrow
\infty$, the area of the critical locus of a special Lagrangian fibration
on $\X_{t_i}$ goes to zero, so its image hopefully deforms to something
of codimension two.

Even once one finds a single special Lagrangian torus, it could fail to
give a fibration either because deformations may not be disjoint from
each other, or the deformations simply may not fill out the entire manifold,
so that $B_i$ has a boundary. The expectation might be that these
sorts of things are more likely to happen near the discriminant locus
in $B_{\infty}$.

(2) We do not expect $B_{0,i}$ to be a Monge-Amp\`ere
manifold, but only an affine K\"ahler manifold. This is because of the
existence of examples of moduli of special Lagrangian tori where this is
not the case: see the work of Matessi in [25].

(3) Stated in the proper way, [11] proves this conjecture for K3 surfaces.

\bigskip

The philosophy of Kontsevich and Soibelman, which I believe is the
right one, is that it may be enough to work purely
with the limiting data. Whereas the original form of the SYZ conjecture
proposed dualizing torus fibrations, we instead dualize the limiting
data. Given a Hessian manifold, one obtains a new affine structure
with local affine coordinates $\check y_i=\partial K/\partial y_i$,
where $K$ is the potential of the metric. The metric remains the
same, but the new potential $\check K$ is the Legendre transform of
$K$. This was first suggested in the context of mirror symmetry
by Hitchin in [16], and this idea was used effectively in [23] and 
[24]\footnote{${}^\dagger$}{Intriguingly, this duality was mentioned in [31],
which gave a reference to a 1985 work in statistics, [2],
which makes serious use of this duality between Hessian manifolds.}.

We are left with two fundamental questions:

\proclaim The Limit Question. How can one calculate the limit data, or
guess it conjecturally either for general degenerations or for standard
cases such as hypersurfaces in toric varieties?

\proclaim The Reconstruction Question. Given a set of limiting data,
how do we reconstruct a family of Calabi-Yau manifolds converging to
this limit?

\medskip
Kontsevich and Soibelman discuss these two questions, suggesting some
approaches involving rigid analytic geometry and Berkovich spaces. From
my point of view, these questions can be developed at a topological level
(where one only preserves the limiting information of monodromy about
$\Delta$), a symplectic level (pay attention only to the affine structure)
and the full metric level. For the limit question, [10] gives a 
conjectural limit for the quintic and its mirror on the topological level,
while [29] gives it for general toric hypersurfaces. The affine
structure can be guessed at from the ideas in [29], and in future work
I will give a purely combinatorial description of a conjectural affine
structure in the limit for hypersurfaces in toric varieties.

For the reconstruction question, the results of [10] allow a 
reconstruction of the underlying topological manifold from the limiting
data in sufficiently generic cases, while current work in progress
of my own explores the
symplectic reconstruction problem.

However, solving these general questions at the metric (or complex structure)
level will require some substantial new ideas.
\bigskip
{\hd Bibliography}

\item{[1]} Altmann, K., ``The Versal Deformation of an Isolated
Toric Gorenstein Singularity,'' {\it Inv. Math.} {\bf 128}, (1997) 443--479.
\item{[2]} Amari, S., {\it Differential-Geometric Methods in Statistics,}
{\it Lecture Notes in Statistics,} {\bf 28}, Springer-Verlag, 1985.
\item{[3]} Audin, M., {\it The Topology of Torus Actions on Symplectic
Manifolds,} {\it Prgress in Mathematics, 93}, 
Birkh\"auser Verlag, Basel, 1991.
\item{[4]} Chiang, T.-M., Klemm, A., Yau, S.-T. and Zaslow, E.,
``Local Mirror Symmetry: Calculations and Interpretations,'' preprint,
hep-th/9903053.
\item{[5]} Givental, A., ``Homological Geometry and Mirror Symmetry,''
{\it Proceedings of the ICM, Zurich, 1994,} (Birkhauser, 1995), 472--480.
\item{[6]} Goldstein, E., ``Calibrated Fibrations,'' preprint,
math/9911093.
\item{[7]} Goldstein, E., ``Calibrated Fibrations on Complete
Manifolds via Torus Action,'' preprint, math/0002097.
\item{[8]} Gross, M., ``Special Lagrangian Fibrations I: Topology,''
in {\it Integrable Systems and Algebraic Geometry},
eds. M.-H. Saito, Y. Shimizu and
K. Ueno, World Scientific, 1998, 156--193.
\item{[9]} Gross, M., ``Special Lagrangian Fibrations II: Geometry,''
{\it Surveys in Differential Geometry,} Somerville: MA, International Press,
1999, 341--403.
\item{[10]} Gross, M., ``Topological Mirror Symmetry,'' preprint, math.AG/9909015
(1999), to appear in {\it Inv. Math}.
\item{[11]} Gross, M., and Wilson, P.M.H., ``Large Complex Structure Limits
of K3 Surfaces,'' preprint, math.DG/0008018.
\item{[12]} Harvey, R., and Lawson, H.B. Jr.,  ``Calibrated Geometries,'' {\it
Acta
Math.} {\bf 148}, 47-157 (1982).
\item{[13]} Haskins, M., ``Special Lagrangian Cones,'' preprint, math.DG/0005164.
\item{[14]} Heinzner, P., and Huckleberry, A., ``K\"ahlerian
Potentials and Convexity Properties of the Moment Map,'' {\it Inv. Math.}
{\bf 126}, (1996) 65--84.
\item{[15]} Hilgert, J., Neeb, K.-H., and Plank, W., ``Symplectic
Convexity Theorems and Coadjoint Orbits,'' {\it Comp. Math.} {\bf 94},
(1994), 129--180.
\item{[16]} Hitchin, N., ``The Moduli Space of Special Lagrangian 
Submanifolds,'' {\it Ann. Scuola Norm. Sup. Pisa Cl. Sci. (4)}, {\bf 25}
(1997) 503--515.
\item{[17]} Hori, K., and Vafa, C., ``Mirror Symmetry,'' preprint, hep-th/0002222.
\item{[18]} Hori, K., Iqbal, A., and Vafa, C., ``$D$-branes and Mirror Symemtry,''
preprint, hep-th/0005247.
\item{[19]} Joyce, D., ``Asymptotically Locally Euclidean Metrics
with Holonomy $SU(m)$,'' preprint, math.AG/9905041.
\item{[20]} Joyce, D., ``Quasi-ALE Metrics with Holonomy $SU(m)$ and
$Sp(m)$,'' preprint, math.AG\-/9905043.
\item{[21]} Joyce, D., ``Special Lagrangian $m$-folds in ${\bf C}^m$ with
Symmetries,'' preprint, math.DG\-/0008021.
\item{[22]} Joyce, D., ``Singularities of Special Lagrangian Fibrations
and the SYZ Conjecture,'' preprint, math.DG/0011179.
\item{[23]} Kontsevich, M., and Soibelman, Y., ``Homological Mirror Symmetry
and Torus Fibrations,'' preprint, math.SG/0011041.
\item{[24]} Leung, N.C., ``Mirror Symmetry Without Corrections,''
preprint, math.DG/0009235.
\item{[25]} Matessi, D., ``Some Families of Special Lagrangian Tori,''
preprint, math.DG/0011061.
\item{[26]} McLean, R.C., `` Deformations of Calibrated Submanifolds,'' 
{\it Comm. Anal. Geom.} {\bf 6}, (1998) 705--747.
\item{[27]} Mikhalkin, G., ``Real Algebraic Curves, the Moment Map,
and Amoebas,'' {\it Ann. of Math.}, {\bf 151}, (2000), 309--326.
\item{[28]} Reid, M., ``Decomposition of Toric Morphisms,'' {\it Arithemetic
and Geometry, Vol. II}, 395--418, {\it Progr. Math. 36}, Birkh\"auser
Boston, 1983.
\item{[29]} Ruan, W.-D., ``Lagrangian Torus Fibration and Mirror Symmetry
of Calabi-Yau Hypersurface in Toric Variety,'' preprint, math.DG/0007028.
\item{[30]} Ruan, W.-D., ``Newton Polygon and String Diagram,'' preprint,
math.DG/0011012.
\item{[31]} Shima, H., and Yagi, K., ``Geometry of Hessian Manifolds,''
{\it Differential Geom. Appl.} {\bf 7}, (1997), 277--290.
\item{[32]} Stenzel, M., ``Ricci-flat Metrics on the Complexification
of a Compact Rank One Symmetric Space,'' {\it Manuscripta Math.} {\bf 80},
91993), 151--163.
\item{[33]} Strominger, A., Yau, S.-T., and Zaslow, E.,  ``Mirror Symmetry is
T-Duality,'' {\it Nucl. Phys.} {\bf B479}, (1996) 243--259.
\item{[34]} Tian, G., and Yau, S.-T., ``Complete K\"ahler manifolds with
zero Ricci curvature, I.''
{\it J. Amer. Math. Soc.} {\bf 3}, (1990) 579--609.
\item{[35]} Tian, G., and Yau, S.-T., ``Complete K\"ahler manifolds with
zero Ricci curvature, II.''
{\it Invent. math.} {\bf 106}, (1991) 27--60.
\item{[36]} Viro, O., ``Gluing of Plane Real Algebraic Curves and Constructions
of Curves of Degrees 6 and 7,'' in {\it Topology (Leningrad, 1982)}, 187--200,
{\it Lecture Notes in Math., 1060}, Springer-Verlag, 1984.
\end